\documentclass[11pt]{amsart}
\usepackage{eucal}
\usepackage{amssymb}
\usepackage{epsfig}
\usepackage{amscd}
\def\P{\mathbb{P}}

\newtheorem{thm}{Theorem}[section]
\newtheorem{Prop}[thm]{Proposition}
\newtheorem{lemma}[thm]{Lemma}

\newtheorem{co}[thm]{Corollary}
\newtheorem{Def}[thm]{Definition}

\newcommand{\Sa}{\Sigma}

\newcommand{\bc}{{\mathbb C}}
\newcommand{\bF}{{\mathbb F}}

\newcommand{\bh}{{\mathbb H}}

\newcommand{\br}{{\mathbb R}}

\newcommand{\ra}{\rightarrow}

\newcommand{\col}{\!:\!}

\newcommand{\fh}{H_{\bF}^}%{\bh_{\bF}^}
\newcommand{\hh}{H_{\bh}^}%{\bh_{\bh}^}
\newcommand{\rh}{H_{\br}^}%{\bh_{\br}^}
\newcommand{\id}{\operatorname{id}}

\newcommand{\im}{\operatorname{Im}}
\newcommand{\re}{\operatorname{Re}}

\newenvironment{pf}{\begin{trivlist}\item[]{\bf Proof:\ }}
{\mbox{}\hfill\rule{.08in}{.08in}\end{trivlist}}

\begin{document}
\title[Deformation of quaternionic  space]
{Deformation space of a non-uniform 3-dimensional real hyperbolic
lattice in quaternionic hyperbolic plane}
\author{Inkang Kim}
\date{}%May 12, 2004}
\maketitle

\begin{abstract}
In this note, we study deformations of a non-uniform real hyperbolic
lattice in quaternionic hyperbolic spaces. Specially we show that
the representations of the fundamental group of the figure eight
knot complement into $PU(2,1)$, cannot be deformed in $PSp(2,1)$ out
of $PU(2,1)$ up to conjugacy.
 \end{abstract}
%  \footnotetext[1]{1991 {\sl{Mathematics Subject Classification.}}51M10, 57S25.}
%    \footnotetext[2]{{\sl{Key words and phrases.}} Quaternionic hyperbolic
%   space, Cayley hyperbolic plane, symmetric rank one spaces,
%     Cartan angular invariant, Toledo invariant, parallel 4- and 8-forms, Gromov bounded
%cohomology, quasifuchsian representations, bendings, rigidity
%\hfil\hfil\hfil}
  \footnotetext[3]{The author
gratefully acknowledges the partial support of KRF grant
(0409-20060066) and a warm support of IHES during his stay.}

\section{Introduction}\label{intro}
In 1960's, A.  Weil  \cite{Weil} proved a local rigidity of a
uniform lattice $\Gamma\subset G$ inside $G$, i.e., he showed that
$H^1(\Gamma,\mathfrak{g})=0$ for any semisimple Lie group $G$ not
locally isomorphic to $SL(2,\br)$. This result implies that the
canonical inclusion map $i:\Gamma \hookrightarrow G$ is locally
rigid up to conjugacy. In other words, for any local deformation
$\rho_t:\Gamma\ra G$ such that $\rho_0=i$, there exists a continuous
family $g_t\in G$ such that $\rho_t= g_t \rho_0 g_t^{-1}$. Weil's
idea is further explored by many others but notably by Raghunathan
\cite{Ra} and Matsushima-Murakami \cite{MM}. Much later Goldman and
Millson \cite{GM} considered the embedding of a uniform lattice
$\Gamma$ of $SU(n,1)$
$$\Gamma \hookrightarrow SU(n,1)\hookrightarrow SU(m,1),\ m>n$$ and
proved that there is still a local rigidity inside $SU(m,1)$ even if
one enlarges the target group. More recently  further examples of
the local rigidity of a complex hyperbolic lattice in quaternionic
K\"ahler manifolds are found in \cite{KKP}
$$\Gamma\hookrightarrow SU(n,1)\subset Sp(n,1)\subset SU(2n,2)
\subset SO(4n,4).$$

But all these examples deal with the standard inclusion map $\Gamma
\hookrightarrow G'$ to use the Weil's original idea about
$L^2$-group cohomology. It has not been much studied yet when
$\rho:\Gamma\ra G'$ is an arbitrary representation which is not an
inclusion.
 In this
note, we study deformations of a non-inclusion representation
$\rho_0$ of a non-uniform lattice $\Gamma$ of $PSL(2,\bc)$ in
semisimple Lie groups $PU(2,1)$ and $PSp(2,1)$
$$\Gamma \stackrel{\rho_0}{\longrightarrow}  PU(2,1) \subset PSp(2,1).$$

This is a sequel to the
 previous paper
\cite{KP} where the deformation of the standard inclusion
$\Gamma\hookrightarrow SO(3,1)\subset SO(4,1)\hookrightarrow
Sp(n,1)$ is studied but techniques are quite different. In
\cite{KP}, we used the group cohomology to prove the local rigidity
following the path of \cite{Weil, Ra}. Here we use explicit
coordinates and matrix calculations to prove a kind of a local
rigidity, namely representations into $PU(2,1)$ cannot be deformed
into $PSp(2,1)$ nontrivially. We calculate the dimension of a
representation variety of the fundamental group of the figure eight
knot complement in complex and quaternionic hyperbolic 2-plane using
Thurston's idea.

In \cite{Falbel}, Falbel constructed a special Zariski dense
discrete representation $\rho_0$ in $PU(2,1)$ with purely parabolic
holonomy for
 a peripheral group. If $g_1,g_2$ are generators of
the fundamental group of the figure eight knot complement, their
images in $SU(2,1)$
\[ G_1=\left [ \begin{matrix}
     1 &   1   & \frac{-1-i\sqrt 3}{2}\\
     0 &   1   & -1  \\
     0 & 0 & 1 \end{matrix} \right]\]
     and
\[ G_2=\left[\begin{matrix}
   1 & 0 & 0 \\
   1 & 1 & 0\\
\frac{-1-i\sqrt 3}{2}&-1& 1 \end{matrix}\right]\] give such a
special representation in $PU(2,1)$. Note that this representation
is not faithful.

More precisely we prove:
\begin{thm}\label{Qstructure}Let $M$ be a figure eight knot complement which can be made up of two
ideal tetrahedra in the quaternionic hyperbolic plane $H^2_\bh$
glued up along faces properly. For the space of representations
$\rho:\pi_1(M)\ra PSp(2,1)$  which do not stabilize a quaternionic
line,
%such that the holonomies around two
%edges multiply to be a pure parabolic element fixing a common point
%of two edges,
the representation variety  in $PSp(2,1)$
 around a discrete representation $\rho_0$, is  of real 3 dimension up to conjugacy.
The representation variety into $PU(2,1)$ around $\rho_0$ is of
dimension 3 up to conjugacy as well. The variety around the
conjugacy class $[\rho_0]$ is parameterized by the three angular
invariants of faces of one ideal tetrahedron.
\end{thm}
%Due to Lemma \ref{pureparabolic}, two pure parabolic elements can
%commute.
%\begin{co}If a (discrete) representation not conjugate into $Sp(1,1)$ sends
%the peripheral subgroup into pure parabolics, then it is  conjugate
%into $SU(2,1)$.
%\end{co}
%\begin{pf}The product of holonomies around two edges fixes the
%common point of two edges, so it is a parabolic element fixing the
%common point. Obviously holonomies of peripheral subgroup fix this
%common point. So the product belongs to the peripheral
%subgroup???(it is true for discrete representations). Then by the
%assumption, it is pure parabolic.
%\end{pf}

 As a corollary we obtain
\begin{co}Any representation from the fundamental group of the
figure eight knot complement into $PU(2,1)$ near the discrete
representation $\rho_0$, cannot be deformed in $PSp(2,1)$ out of
$PU(2,1)$ up to conjugacy.
\end{co}

{\bf Acknowledgements} The author thanks an anonymous  referee for
correcting some errors of the first version of the paper and
constructive suggestions.

\section{Preliminaries}\label{Pre}
\subsection{Different models of hyperbolic space}
The set $\bh$ of quaternions are  $\{x=x_1+ix_2+jx_3+kx_4|
x_i\in\br\}$ with the multiplication law $ij=k,\ jk=i,\
i^2=j^2=k^2=-1$. We set $\Im x=ix_2+jx_3+kx_4$ and $\bar
x=x_1-ix_2-jx_3-kx_4$. We call $x$ pure imaginary  if $\bar x=-x$.
Quaternion number is a non-commutative division ring and by the
abuse of notations, we will set $x^{-1}=\frac{1}{x}$ to be the
multiplicative inverse of $x$. Up to section \ref{qstructure},
multiplication by quaternions on $\bh^n$ is {\bf on the left} and
matrices act {\bf on the right}. Let $J_0$ be
\[ \left [ \begin{matrix}
         0   &  0 & 1 \\
         0   & I_{n-1}  & 0 \\
         1   & 0   & 0  \end{matrix} \right ]. \]
Define $\langle Z, W \rangle_0= Z J_0 W^*$ where
$Z=(z_1,\cdots,z_{n+1})\in \bh^{n+1}$.
 Then $A\in Sp(n,1)$ if
$$AJ_0A^*=J_0.$$
Hence $$A^{-1}=J_0A^*J_0.$$

 One can define projective
models, called {\bf Siegel domains}, of the hyperbolic spaces $\fh
n,\ \bF=\bc,\bh$ as the set of negative lines in the Hermitian
vector space $\bF ^{n,1}$, with Hermitian structure given by the
indefinite $(n,1)$-form
$$
\langle Z, W\rangle_0= Z J_0 W^*.
$$  Namely $\fh n$ is the left projectivization $\P V_-$ of the set
$$V_-=\{Z\in \bF^{n,1}: \langle Z, Z \rangle_0 <0\}.$$
The boundary of the Siegel domain consists of projectivized zero
vectors $$V_0=\{Z\in \bF^{n,1}-\{0\}: \langle Z, Z \rangle_0=0\}$$
together with a distinguished point at infinity $\infty$. The finite
points in the boundary carry the structure of the generalized
Heisenberg group $\bF^{n-1}\times \Im \bF$ with the group law
$$(Z,t)(W,s)=(Z+W, t+s- 2\Im \langle\langle Z, W \rangle\rangle )$$
where $\langle\langle Z, W \rangle\rangle=ZW^*=\sum z_i\bar w_i$ is
the standard positive definite Hermitian form on $\bF^{n-1}$.
Motivated by this one can define {\bf horospherical coordinates} for
$H^n_\bF$
$$\{(z,t,v)\in \bF^{n-1}\times \Im \bF\times \br_+\}.$$

From now on we will take $n=2$ so that we will deal with only two
dimensional hyperbolic spaces. A coordinate change $\psi$ from
Heisenberg coordinates $(z,t),z\in \bF, t\in \Im \bF$ to the
boundary of Siegel domain is
\begin{eqnarray}\label{coo}
(\frac{-|z|^2+t}{2},z,1)\end{eqnarray} with one extra equation
$$\psi(\infty)=(1,0,0).$$

If $U\in Sp(1), \mu\in Sp(1), r\in\br^+$, the action fixing $0$ and
$\infty$ is given by $$ (z,t)\ra (r\mu^{-1} z\mu U, r^2\mu^{-1} t
\mu).$$ See \cite{Kim,KParker}. In matrix form acting on the right
\[ H_{\mu,U,r}=\left [ \begin{matrix}
         r\mu  &  0 & 0 \\
         0   & \mu U & 0 \\
         0   & 0   & \frac{\mu}{r}  \end{matrix} \right ] \]
So the hyperbolic isometry fixing $\infty$ and $0$ is determined by
$\mu,\nu=\mu U\in Sp(1)$ and $r\in\br^+$, so it is 7 dimensional.
For complex hyperbolic space $H^2_\bc$, $\mu=1$ and $U\in U(1)$.
\begin{lemma}\label{fix}The set of isometries fixing three points on the ideal
boundary of $H^2_\bh$, which do not lie on the quaternionic line, is
one dimensional whereas it is unique in $H^2_\bc$.
\end{lemma}
\begin{pf}We may assume that three points are
$$\infty,0,(1,it)$$ up to the action of $Sp(2,1)$.
If $H_{\mu,U,r}$ fixes $(1,it)$, then $U=1,r=1$ and
$\mu^{-1}it\mu=it$. It is easy to show that $\mu=e^{i\theta}$,
showing that $\mu$ has one degree of freedom.
\end{pf}

 The Heisenberg group acts
by right multiplication:
$$T_{(z,t)}(\zeta,v)=(z+\zeta, t+v- 2 \Im \zeta\bar z).$$
 In matrix form acting on $\bh^{2,1}$ on the right
\[T_{(z,t)}= \left [ \begin{matrix}
         1  &  0 & 0 \\
         -\bar z   & 1  & 0 \\
         \frac{-|z|^2+t}{2}   & z   & 1  \end{matrix} \right ] \]

Then a hyperbolic isometry fixing $\infty$ and $(z,t)$ is
\[ T_{(-z,-t)}\circ H_{\mu,U,r}\circ T_{(z,t)}=\left[
\begin{matrix}
   r\mu &0 &0 \\
   r\bar z \mu-\nu\bar z &\nu & 0 \\
   r\frac{-|z|^2-t}{2}\mu+z\nu\bar z+\frac{\mu}{r}\frac{-|z|^2+t}{2}&-z\nu+\frac{\mu}{r}z&\frac{\mu}{r}
    \end{matrix}\right]\]
where $\nu=\mu U\in Sp(1)$. This group is also 7 dimensional
determined by $\mu,\nu\in Sp(1),r\in\br^+$. We call $T_{(z,t)}$ a
pure parabolic whereas $H_{(\mu,U,1)}\circ T_{(z,t)}$
ellipto-parabolic if it fixes a unique point at infinity.
\begin{lemma}\label{pureparabolic}Two pure parabolic elements $T_{(z,t)},T_{(w,s)}$
commute if $w\bar z$ is real, i.e. $w=rz$ for some real $r$.
\end{lemma}
\begin{pf}A direct calculation shows that two elements commute iff
$w\bar z=z\bar w$, which implies that $w\bar z$ is real.
\end{pf}
 There is one more isometry interchanging $\infty$ and $0$
whose matrix form is
\[ \left [ \begin{matrix}
   0 & 0 & 1 \\
   0 & 1 & 0 \\
   1 & 0 & 0 \end{matrix} \right] \] and in Heisenberg coordinates
$$I(z,t)=(-\frac{2}{|z|^2+\bar t}z, 4\frac{\bar{t}}{|z|^4+|t|^2}).$$
          A reflection with respect to $H^2_\br$ in $H^2_\bh$ is
an isometry. In Heisenberg coordinates, it is
 $$(z,t)\ra (\bar z, \bar t).$$
 So
 $$(z,t)\ra (\frac{-2}{|z|^2+{ t}}\bar z,\frac{{4t}}{|z|^4+|t|^2})$$
 is an isometry interchanging $\infty$ and $0$.
Composing with $ (z,t)\ra (r\mu z U\mu^{-1}, r^2 \mu t \mu^{-1})$ we
get
$$R(z,t)=(r\mu\frac{-2 }{|z|^2+{ t}}\bar z\nu,r^2\mu \frac{{4t}}{|z|^4+|t|^2}\mu^{-1}),$$
where $r\in \br^+, \mu,\nu\in Sp(1)$.

%Up to $PSp(n,1)$, four-tuples of points on the ideal boundary of quaternionic hyperbolic space
%are uniquely determined by their Cartan angular invariants.

\subsection{Angular invariant}
To define angular invariant we introduce the {\bf unit ball model}
$\{Z\in \bF^n:||Z||<1\}$ to make it compatible with existing
literatures where $\bF^n$ is equipped with the standard positive
definite Hermitian form. Two points $(0',-1)$ and $(0',1)$ will play
a special role. There is a natural map from a unit ball model to
$\P(\bF^{n,1})$ where $\bF^{n,1}$ is equipped with a standard
$(n,1)$ Hermitian product
$$ \langle z,w\rangle =z_1\overline
w_1+\cdots+z_n\overline w_n-z_{n+1}\overline w_{n+1}\,,
$$ defined as
$$(w',w_n)\ra (w',w_n,1).$$ From unit ball model to the horospherical model, one defines the coordinates change as
$$(z',z_n)\ra (\frac{z'}{1+z_n},\frac{2\Im z_n}{|1+z_n|^2},\frac{1-|z_n|^2-|z'|^2}{|1+z_n|^2}).$$
Its  inverse from the horospherical model  to $\P\bF^{n,1}$ is given
by
\begin{eqnarray}\label{ch}(\xi,v,u)=[(\xi,\frac{1-|\xi|^2-u+v}{2},\frac{1+|\xi|^2+u-v}{2})],\end{eqnarray}
where $v$ is pure imaginary, i.e., $iv$ in complex case, and
$iv_i+jv_2+kv_3$ in quaternionic case. Note this coordinate change
is different from (\ref{coo}) since we used a different Hermitian
product. According to this coordinate change, $(0',1)=[(0',1,1)]$
corresponds to the identity element $(0',0)$ in Heisenberg group,
$(0',-1)=[(0',-1,1)]$ to $\infty$, and $(0',0)$ to $(0,0,1)$.
\begin{Def}\label{cart} The complex Cartan angular invariant $\mathbb A(x_1,x_2,x_3)$ of the ordered triples $(x_1,x_2,x_3)$
in $\partial H^n_\bc$ is introduced by Cartan \cite{Ca} and defined
to be the argument between $-\frac{\pi}{2}$ and $\frac{\pi}{2}$ of
the Hermitian triple product
$$ -\langle\tilde{x_1},\tilde{x_2},\tilde{x_3} \rangle= -\langle
\tilde{x_1},\tilde{x_2}\rangle \langle \tilde{x_2},\tilde{x_3}
\rangle \langle\tilde{x_3},\tilde{x_1}\rangle\in\bc\ $$  where
$\tilde x_i$ is a lift of $x_i$ to $\bc^{n,1}$. It can be obtained
by, up to constant, integrating the K\"ahler from on $H^n_\bc$ over
the geodesic triangle spanned by three points \cite{Do}, hence it is
a bounded cocyle. It satisfies the cocycle relation: for
$(x_1,x_2,x_3,x_4)\in\partial H^n_\bc$
\begin{eqnarray}\label{cocycle}
\mathbb A(x_1,x_2,x_3) +\mathbb A(x_1,x_3,x_4)  =\mathbb
A(x_1,x_2,x_4)+\mathbb A(x_2,x_3,x_4).
\end{eqnarray}

 The quaternionic Cartan angular invariant of a triple $x$,
$0\leq {\mathbb A}_{\bh}(x)\leq \pi/2$,  is the angle between the
first coordinate line $\br e_1=(\br,0,0,0)\subset\br^4$ and the
Hermitian triple product
$$
\langle\tilde{x_1},\tilde{x_2},\tilde{x_3} \rangle= \langle
\tilde{x_1},\tilde{x_2}\rangle \langle \tilde{x_2},\tilde{x_3}
\rangle \langle\tilde{x_3},\tilde{x_1}\rangle\in\bh\,,
$$
where we identify ${\mathbb H}$ and ${\mathbb R}^4$.
\end{Def}
Note that the invariant is unchanged under the homothety by nonzero
real numbers, i.e., the triples $x$ and $rx$ have the same angular
invariant.
\begin{Prop}\label{ang}
{\em ( \cite{Pernas},\cite{KimJKMS}).} Let $x=(x_1,x_2,x_3)$ and
$y=(y_1,y_2,y_3)$ be pairs of distinct triples of points in $\hh n$.
Then $ {\mathbb A}_{\bh}(x)={\mathbb A}_{\bh}(y) $ if and only if
there is an isometry $f\in PSp(n,1)$ such that $f(x_i)=y_i$ for
$i=1,2,3$.
\end{Prop}

\begin{pf}
Applying a homothety by nonzero real number, we may assume that our
triples have Hermitian products $X=\langle
\tilde{x_1},\tilde{x_2},\tilde{x_3} \rangle$ and $ Y=\langle
\tilde{y_1},\tilde{y_2},\tilde{y_3}\rangle$ with $|X|=|Y|$.

Now let us assume that ${\mathbb A}_{\bh}(x)={\mathbb A}_{\bh}(y)$.
This and $|X|=|Y|$ imply that there is an orthogonal transformation
$M \in SO(3)\times \{\id\}$ acting on $\bh=\br^4$ that leaves
invariant the real axis in $\bh$ and maps $X$ to $Y$.  Since the
conjugation action of $Sp(1)$ in $\bh$ is $SO(3)$ action, there is
$\mu \in Sp(1)$ such that
$$
\langle \tilde{x_1},\tilde{x_2},\tilde{x_3} \rangle= \mu \langle
\tilde{y_1},\tilde{y_2},\tilde{y_3} \rangle \bar{\mu}\,.
$$
To finish the proof it is enough to choose lifts $\tilde{x_i}$ and
$\tilde{y_i}$ of points $x_i$ and $y_i$, $i=1,2,3$, so that $\langle
\tilde{x_i},\tilde{x_j}\rangle=
 \langle \tilde{y_i},\tilde{y_j}\rangle$. Indeed,
 then there is $A\in Sp(n,1)$ such that $A(\tilde{x_i})=\tilde{y_i}$, $i=1,2,3$.
 Then it descends to an element $f\in PSp(n,1)$ such that
$f(x_i)=y_i$ for $i=1,2,3$.

To obtain those lifts, we first replace $\tilde{y_1}$ by $\mu
\tilde{y_1}$ (still denote it by
 $\tilde{y_1}$) and
 get $\langle \tilde{x_1},\tilde{x_2}\rangle
 \langle \tilde{x_2},\tilde{x_3} \rangle \langle\tilde{x_3},\tilde{x_1}\rangle=
 \langle \tilde{y_1},\tilde{y_2}\rangle
 \langle \tilde{y_2},\tilde{y_3} \rangle \langle\tilde{y_3},\tilde{y_1}\rangle$.
  Replacing $\tilde{x_2}$ and $\tilde{x_3}$ by $\mu_2\tilde{x_2}$ and $ \mu_3\tilde{x_3}$
   if necessary, we can make
   $\langle\tilde{x_2}, \tilde{x_3}\rangle=\langle \tilde{y_2}, \tilde{y_3}
\rangle$ and $ \langle \tilde{x_3},\tilde{x_1}\rangle=\langle
\tilde{y_3},\tilde{y_1}
   \rangle$. Now the equation becomes
   $$
\langle \tilde{x_1},\tilde{x_2}\rangle
 \langle \tilde{x_2},\tilde{x_3} \rangle \langle\tilde{x_3},\tilde{x_1}\rangle=
    |\mu_2|^2|\mu_3|^2 \langle \tilde{y_1},\tilde{y_2}\rangle
\langle \tilde{y_2},\tilde{y_3}\rangle
\langle\tilde{y_3},\tilde{y_1}\rangle\,,
$$
and we get $\langle \tilde{x_1},\tilde{x_2}\rangle= r\langle
\tilde{y_1},\tilde{y_2}\rangle$ where $r=|\mu_2||\mu_3|$. Then
replacing $\tilde{x_1}$, $\tilde{x_2}$, $\tilde{x_3}$ and
$\tilde{y_1}$ by $r^{-1}\tilde{x_1}$, $r^{-1}\tilde{x_2}$,
$r\tilde{x_3}$ and $r^2\tilde{y_1}$ respectively, we finally get
$\langle \tilde{x_i},\tilde{x_j}\rangle= \langle
\tilde{y_i},\tilde{y_j}\rangle$, and hence a desired $f\in
PSp(n,1)$.

The converse is trivial.
\end{pf}

\begin{thm}\label{dist}
For distinct points $x_1,x_2,x_3 \in \partial{H^n_{\mathbb H}}$, let
$\sigma_{12}$ and $\Sigma_{12}$ be real and quaternionic geodesics
containing the two points $x_1$ and $x_2$, and $\Pi\col H^n_{\mathbb
H}\rightarrow \Sigma_{12}$ be the orthogonal projection. Then
$$
|\tan{{\mathbb A}_{\bh}(x)}|=\sinh(d(\Pi{x_3},\sigma_{12}))
$$
where $d$ is the hyperbolic distance in $H^n_{\mathbb H}$.
\end{thm}
\begin{pf}
        Up to an isometry (in the unit ball model of $\hh n$), we may assume that
the triple $x$ consists of $ x_1=(0,-1)$, $x_2=(0,1)$, $x_3=(z',z_n)
$,\,\, whose lifts are $ \tilde{x_1}=(0,-1,1), \tilde{x_2}=(0,1,1),
\tilde{x_3}=(z',z_n,1) $.
  In this setting $\sigma_{12}=\{(0,t)\col t\in{\mathbb R},|t|<1\},
\Sigma_{12}=\{(0,z)\col z\in{\mathbb H},|z|<1\}$, and $\langle
\tilde{x_1},\tilde{x_2},\tilde{x_3} \rangle=2(\bar{z_n}-1)(1+z_n)$.
So we get
$$
\big|\tan{{\mathbb A}_{\bh}(x)}\big|=\frac{|2\im
(z_n)|}{1-|z_n|^2}\,.
$$

On the other hand, we note that $\Pi(x_3)=z_n$ and $\Sigma_{12}$ has
the Poincar\'e ball model geometry of $\rh 4$ with sectional
curvature $-1$. Choose a hyperbolic two plane in $\Sa_{12}$ that
contains the geodesic $\sigma_{12}$ and $z_n$. This plane is a
Poincar\'e disk with curvature $-1$, where we can write
$z_n=\re{z_n}+i|\im{z_n}|$.   Let $d$ be the hyperbolic distance
between the point $z_n$ and the real axis in that Poincar\'e disk.
Then a direct calculation shows that
$\sinh(d)=|2\im(z_n)|/(1-|z_n|^2)$.
\end{pf}

\section{Q-structure}\label{qstructure}
After the complete hyperbolic structure on the figure eight knot
complement was first given in \cite{Ri}, W. Thurston describes a
complete, finite volume hyperbolic structure on the complement of
the figure eight knot in the 3-sphere  in \cite{Th} by gluing two
tetrahedra. Then, he shows how to deform it to non-complete
structures. All these structures have holonomies, which are
homomorphisms of the fundamental group $\Gamma$ of the figure eight
knot complement (a 2 generators and 1 relator group) to $SO(3,1)^0$.
In fact, the {\bf character variety}
$\chi(\Gamma,SO(3,1)^0)=Hom(\Gamma,SO(3,1)^0)/SO(3,1)^0$ is a smooth
1-dimensional complex manifold near the conjugacy class of the
holonomy of the complete hyperbolic structure.

Let $\rho =\Gamma\ra PSp(2,1)$ be the holonomy of the complete
hyperbolic structure, followed by the embedding $SO(3,1)^0 \ra
SO(4,1)^0 =PSp(1,1)\ra PSp(2,1)$. In \cite{KP}, it is shown that any
parabolicity-preserving local deformation around $\rho$ preserves a
quaternionic line again. It is conjectured that any deformation
around $\rho_0$ preserves a quaternionic line.
%We show that the character variety
%$\chi(\Gamma,PSp(2,1))$ has positive dimension near the conjugacy
%class of $\rho_0$, and that most deformations of $\rho_0$ are
%Zariski dense.

In this note, we prove that the component, containing $\rho_0$
(introduced in section \ref{intro}) with purely parabolic holonomy
for
 a peripheral group, of
the space of representations, which do not stabilize a quaternionic
line, of the fundamental group of the figure eight knot complement
in $PSp(2,1)$ is actually conjugate into $PU(2,1)$.
\subsection{The figure eight knot complement}

Consider the following complex. Glue together two tetrahedra by
identifying faces pairwise according to the pattern indicated in
Figure 1.
% (there is a unique set of identifications which respect the
%arrows, it maps face $A$ to $A'$, $B$ to $C'$, $C$ to $B'$ and $D$
%to $D'$, reversing orientation on faces).
%\begin{figure}
%\includegraphics[width=.5\linewidth]{figureeight.pdf}
%\end{figure}

The resulting complex has two 3-cells, four faces, two edges and one
vertex. According to W. Thurston, the complement of the vertex is
homeomorphic to the complement $M$ of the figure eight knot in the
3-sphere. Identify the two 3-cells of $M$ with regular ideal
tetrahedra in (compactified) $H_{\br}^{3}$. This defines a
hyperbolic structure on $M$, and therefore a homomorphism $\rho_0
:\Gamma=\pi_1 (M)\ra SO(3,1)^0$. Here is how W. Thurston deforms it.
Up to isometry, an ideal tetrahedron in $H_{\br}^{3}$ is
characterized by one complex number. Identifying the 3-cells of $M$
with arbitrary ideal hyperbolic tetrahedra defines a hyperbolic
structure on the complement of the 2 skeleton, depending on two
complex parameters. One has to make sure that the gluing maps are
isometries which extend the hyperbolic structure across the faces.
In order for the hyperbolic structure to extend across the two
edges, two algebraic equations must be satisfied, but one of them
turns out to follow from the other as we can see from Lemma
\ref{pure}. As a consequence, one obtains a (complex) one parameter
family of hyperbolic structures on $M$.

We adapt this construction to obtain homomorphisms $\Gamma\ra
PSp(2,1)$. For this, we first classify ideal tetrahedra in
$H_{\bh}^{2}$, then introduce the relevant geometric structure,
baptised $Q$-structure, and describe the compatibility equations
along edges.

\subsection{Ideal triangles and tetrahedra in $H_{\bh}^2$}

The group $PSp(2,1)$ is not transitive on triples of points of
$\partial H_{\bh}^2$. By Proposition \ref{ang}, a pair of triples
are mapped to each other by an isometry iff they have the same
angular invariant.
% Indeed, a pair of points defines a real
%geodesic $\ell$, and a quaternionic line $\ell_{\bh}$. Denote by $z$
%the orthogonal projection of a third point on $\ell_{\bh}$. The
%distance $d(z,\ell)$ is an invariant of the triple. According to L.
%Pernas, \cite{Pernas}, this is the only one. It can be expressed in
%terms of Pernas' {\em angular invariant}, a quaternionic version of
%Cartan's angular invariant in complex hyperbolic geometry. First,
%Pernas defines the {\em argument} of a quaternion $q$. This is the
%real number $\mathrm{arg}(q)\in[0,\pi]$ such that there exists a
%unit imaginary quaternion $i$ such that $q=e^{i\mathrm{arg}(q)}$.
%This is a complete conjugacy invariant of $q$. Let $X_1$, $X_2$ and
%$X_3$ be nonzero null vectors in $\bh^{2,1}$ representing three
%distinct points $x_1$, $x_2$ and $x_3 \in \partial H_{\bh}^2$. Then
%the angular invariant is
%\begin{eqnarray*}
%C(x_1,x_2,x_3)=\mathrm{arg}(-\langle X_1 ,X_2 \rangle\langle X_2
%,X_3 \rangle\langle X_3 ,X_1 \rangle).
%\end{eqnarray*}
%It does not depend on the choice of vectors $X_j$, nor on the order
%of $x_1$, $x_2$ and $x_3$.
By Theorem \ref{dist} the angular invariant of a triple vanishes if
and only if all points sit in a (compactified) totally real totally
geodesic plane. It takes value $\pi/2$ if and only if all points sit
in a (compactified) quaternionic line.

If three ideal points $x_1$, $x_2$ and $x_3 \in \partial H_{\bh}^2$
do not belong to a quaternionic line, the set of isometries that fix
them is one dimensional by Lemma \ref{fix}. It follows that if two
triangles have equal angular invariants different from $\pi/2$,
there is a one dimensional family of isometries that sends one to
the other. Since the boundary of $H_{\bh}^2$ is 7-dimensional, for
each $c\in[0,\pi/2)$, the space of ideal tetrahedra
$(x_1,\ldots,x_4)$ with ${\mathbb A}_{\bh}(x_1,x_2,x_3)=c$ up to
isometry has dimension 6. It follows that ideal tetrahedra up to
isometry depend on 7 parameters.

\begin{Def}Let $A=\{x_1,\cdots,x_k\},\ k\geq 3,$ be a disjoint collection of
ideal points on the ideal boundary of a rank one symmetric space
$X$. The geometric center of $A$ in $X$ is the barycenter of the
associated measure $\delta_A=\sum \delta_{x_i}$. In more details,
let
$$F(x)=\int_{\partial X} B_o(x,\xi)d\delta_A(\xi)$$ be a function
defined on $X$ where $B_o$ is the Busemann function normalized that
$B_o(o,\xi)=0$. Then it is strictly convex and its value goes to
$\infty$ as $x$ tends to $\partial X$. The barycenter $x_0$ of
$\delta_A$ can be written as
$$dF_{(x_0)}(\cdot)=\int_{\partial X} (dB_0)_{(x_0,\xi)}(\cdot) d\delta_A(\xi)=0.$$
\end{Def}
Let $\Delta$ denote a fixed regular ideal tetrahedron in
$H_{\br}^3$, let $\dot{\Delta}\subset \Delta$ be the complement of
the 1-skeleton. Given an ideal tetrahedron (i.e. 4 distinct points
at infinity $(x_1,\ldots,x_4)$) in $H_{\bh}^2$, the {\em straight
singular simplex} spanning them is the continuous map of
$\dot{\Delta}$ to $H_{\bh}^2$ defined as follows. For each face
$(s_i,s_j,s_k)$ of $\Delta$, map the barycenter $s_{ijk}$ of
$(s_i,s_j,s_k)$ in $H_{\br}^3$ to the geometric barycenter $x_{ijk}$
of $(x_i,x_j,x_k)$ in $H_{\bh}^2$. Map the orthogonal projection of
$s_{ijk}$ to the edge $[s_i,s_k]$ to the orthogonal projection of
$x_{ijk}$ to the geodesic $[x_i,x_j]$ defined by $x_i$ and $x_j$,
extend to an isometric map of edge $[s_i,s_k]$ onto geodesic
$[x_i,x_j]$. Then map each geodesic segment joining $s_{ijk}$ to
$[s_i,s_j]$ to a constant speed geodesic segment joining $x_{ijk}$
to the corresponding point $[x_i,x_j]$. Finally, map each geodesic
segment joining the barycenter of $(s_1,\ldots,s_4)$ to a point on a
face to a constant speed geodesic segment from the barycenter of
$(x_1,\ldots,x_4)$ to the corresponding point in the previously
defined parametrizations of faces. The obtained map being in general
discontinuous along edges, but let us ignore edges.
\subsection{$Q$-structures}

\begin{Def}
\label{Q} Let $M$ be a manifold. A {\em $Q$-structure} on $M$ is an
atlas of charts $\phi_j =U_j \to H_{\bh}^2$ which are continuous
maps from open sets of $M$ to $H_{\bh}^2$, such that on $U_j \cap
U_k$, $\phi_k =\psi_{jk}\circ\phi_j$ for some unique $\psi_{jk}\in
Sp(2,1)$.
\end{Def}

Pick a pair of ideal tetrahedra whose faces have pairwise equal
angular invariants, all  different from $\pi/2$. Map the 3-cells of
the figure eight knot complement $M$ to $H_{\bh}^2$ using the
straight singular simplices spanning chosen ideal tetrahedra. This
defines a $Q$-structure on $M$ with 2-skeleton deleted.

The $Q$-structure extends across faces. Indeed, each face
$(x_i,x_j,x_k)$ of tetrahedron $T$ is isometric to a unique face
$(y_{i'},y_{j'},y_{k'})$ of tetrahedron $T'$. Let $\psi\in Sp(2,1)$
be the unique isometry which maps one face to the other. Then the
straight singular simplices spanning $T$ and $\psi^{-1}(T')$ take
the same values along the common face, and so form a chart defined
in a neighborhood of that face. The elements of $Sp(2,1)$ realizing
the change of charts with the previously defined two charts are
identity and $\psi$ respectively.

The $Q$-structure extends across an edge $[s_i,s_j]$ if and only if
its holonomy, an element of $Sp(2,1)$, around that edge, equals the
identity. Let us compute holonomy based at a point of $T$. A priori,
we know that holonomy maps the image geodesic $[x_i,x_j]$ to itself.
As in section \ref{Pre}, stabilizer of $[x_i,x_j]$ is $\br\times
Sp(1)Sp(1)$. Therefore, vanishing of holonomy amounts to 7
equations. Since there are two edges, we get 14 equations. But we
expect equations provided by the two edges  to be dependent, as it
happens in $SO(3,1)$.

\begin{figure}
$$\includegraphics[width=.8\linewidth]{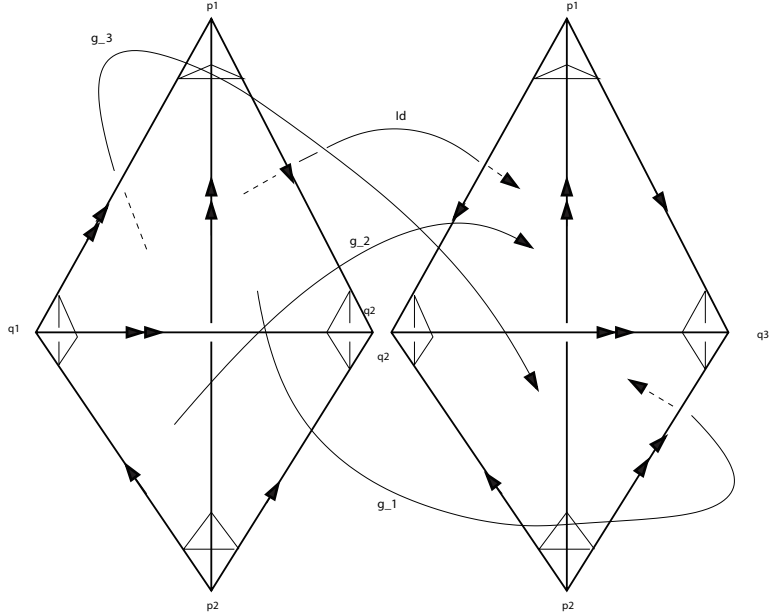}$$
 \caption{Gluing pattern of figure
eight knot complement}
\end{figure}
%\begin{figure}
%$$\includegraphics[scale=.7]{limitset.eps2}$$
%$$\includegraphics[width=.5\linewidth]{limitset2.eps}$$
%\caption{The end points of the word beginning by $x^{-1}, x$ lie in
%the distinct grey regions.}
%\end{figure}
\begin{Prop}Two holonomies around two edges in the complement of
figure eight knot complement obtained by gluing two ideal tetrahedra
are $g_1^{-1}g_3g_2^{-1}g_1g_3^{-1}$ and $g_2^{-1}g_3g_2g_1^{-1}$
where $g_1,g_2,g_3$ are elements in $Sp(2,1)$ appearing in gluing
pattern. One can permute the order of elements appearing in the
products of $g_1,g_2,g_3$.
\end{Prop}
\begin{pf}
By the gluing pattern, we obtain two pictures around two edges and
the proof follows. See Figure 1.
\end{pf}
\begin{lemma}\label{pure}If the holonomies around two edges multiply to be a pure parabolic element
fixing a common point of two edges, which is automatically satisfied
if the holonomies around edges are trivial, then one holonomy is
determined by the other.
\end{lemma}
\begin{pf}Two tetrahedra glue up together to produce two edges $p_1p_2$
and $p_1q_2$. Once we fix two tetrahedra, two holonomy $H_1$ and
$H_2$ around $p_1p_2$ and $p_1q_2$ are hyperbolic isometries
stabilizing them. But since two edges share $\infty$, product
$H_1H_2$ should fix $\infty$. Now put a restriction that $H_1H_2$ is
a pure parabolic isometry. Then
\[ H_1=\left [ \begin{matrix}
         q\alpha &  0 & 0 \\
         0   & \beta & 0 \\
         0   & 0   & \frac{\alpha}{q}  \end{matrix} \right ] \]
         and
\[ H_2=\left[
\begin{matrix}
   r\mu &0 &0 \\
   r\bar z \mu-\nu\bar z &\nu & 0 \\
   r\frac{-|z|^2-t}{2}\mu+z\nu\bar z+\frac{\mu}{r}\frac{-|z|^2+t}{2}&-z\nu+\frac{\mu}{r}z&\frac{\mu}{r}
    \end{matrix}\right]\]
where $q,r\in\br^+, \alpha,\beta,\mu,\nu\in Sp(1)$. Then since
$H_1H_2$ is a pure parabolic element fixing $\infty$,
\[H_1H_2=T_{(w,s)}= \left [
\begin{matrix}
         1  &  0 & 0 \\
         -\bar w   & 1  & 0 \\
         \frac{-|w|^2+s}{2}   & w   & 1  \end{matrix} \right ] \]
for some $(w,s)$. From this equation we obtain
$$qr\alpha\mu=1,\beta\nu=1, \frac{\alpha\mu}{qr}=1.$$
So $r=\frac{1}{q}, \mu=\frac{1}{\alpha},\nu=\frac{1}{\beta}$. This
shows that $H_2$ is completely determined by $H_1$.
\end{pf}
\begin{co}The dimension of the space of real hyperbolic structures
near the complete one on the figure eight knot complement is 2.
\end{co}
\begin{pf}
 An ideal
tetrahedron in a quaternionic line (which is isometric to $H^4_\br$)
is determined by one complex variable. So two tetrahedra have two
complex parameters. A holonomy around an edge belongs to the
stabilizer of a geodesic, in our notation, $H_{\mu,I,r}$ so that
$\mu\in U(1)\subset Sp(1)$. So two holonomies around  edges fix
$\infty$ in common, so its product is naturally a parabolic element.
Then by Lemma \ref{pure}, one holonomy determines the other. Hence
there are two complex parameters with one complex equation and the
solution space is of one complex dimension.
\end{pf}

 First we calculate the dimension of the representation variety near $\rho_0$ in $PU(2,1)$.
\begin{Prop}\label{complex}The dimension of the component of the representation
variety containing $\rho_0$ from the fundamental group of the figure
8 knot complement to $PU(2,1)$ is 3 up to conjugacy.
\end{Prop}
\begin{pf}We claim that to choose an ideal tetrahedron there is 4 degrees of freedom. Choosing three points up to
the action of $PU(2,1)$ is one degree of freedom corresponding to
the Cartan angular invariant. Once three points are fixed, there are
3 degrees of freedom for the last vertex since the boundary of
$H^2_\bc$ is three dimensional. Hence there are total 4 degrees of
freedom to determine an ideal tetrahedron. To determine the second
one, we claim that there is only one degree of freedom. By gluing
pattern, three vertices of the second tetrahedron is determined
according to the angular invariant. The last vertex of the second
tetrahedron is connected to these three vertices to form 3 faces
whose angular invariants are pre-determined by the gluing pattern.
Since one angular invariant is determined if the other three are
known in a tetrahedron by cocycle relation (\ref{cocycle}), two more
angular invariants will determine all the angular invariant. Since
the last vertex can move around 3-dimensional space $\partial
H^2_\bc$ with pre-determined two angular invariants there is only
$3-2$ degree of freedom to choose the second tetrahedron. So there
are total 4+1 degrees of freedom to choose two tetrahedra to glue
them according to the pattern. Then 5 points of two tetrahedra can
be written as
$$p_1=\infty,p_2=0, q_1=(1,t),q_2=(z,s),q_3=(w,r)$$ where $z,w\in\bc, t,s,r\in
\Im \bc$.  A coordinate change from horospherical coordinates
$(z,t),z\in \bc, t\in \Im \bc$ to $\bc^{2,1}$ is
$$(\frac{-|z|^2+t}{2},z,1).$$ Then
$$p_1=(1,0,0),p_2=(0,0,1),q_1=(\frac{-1+t}{2},1,1),q_2=(\frac{-|z|^2+s}{2},z,1),$$$$q_3=(\frac{-|w|^2+r}{2},w,1).$$
As above since there are 5 parameters, there is only one
degree of  freedom for $(w,r)$. This can be easily seen as follows.
Note that the isometries in $PU(2,1)$
$$g_1:(q_2,q_1,p_1)\ra (q_3,p_2,p_1)$$
$$g_2:(p_2,q_1,q_2)\ra (p_1,q_2,q_3)$$
$$g_3:(q_1,p_2,p_1)\ra (q_2,p_2,q_3)$$ are all uniquely determined
by their angular invariants. Hence  $g_1$ gives rise to one  real
equation in $t,z,s,w,r$ which can be derived from the angular
invariant of the faces $(q_2,q_1,p_1)$ and $(q_3,p_2,p_1)$. Indeed
one can do  explicit calculations. Using the coordinate change
formula in (\ref{ch}),
$$p_1=(0,-1,1), p_2=(0,1,1), q_1=(1,\frac{t}{2},\frac{2-t}{2}),$$$$
q_2=(z,\frac{1-|z|^2+s}{2},\frac{1+|z|^2-s}{2}),
q_3=(w,\frac{1-|w|^2+r}{2},\frac{1+|w|^2-r}{2})$$ in $\bc^{2,1}$
with the standard $(2,1)$ Hermitian form $\langle Z, W \rangle=
z_1\bar w_1+ z_2\bar w_2-z_3\bar w_3$ so that
$$-\langle q_2, q_1\rangle\langle q_1, p_1\rangle\langle p_1,q_2\rangle=\frac{|z-1|^2-s-z+\bar z+t}{2}$$
$$-\langle q_3,p_2\rangle\langle p_2,p_1\rangle \langle p_1, q_3 \rangle=|w|^2-r.$$
From $\mathbb A(q_2,q_1,p_1)=\mathbb A (q_3,p_2,p_1)$, we get
\begin{eqnarray}\label{ang}\frac{r}{|w|^2}=\frac{s+z-\bar z-t}{|z-1|^2}.\end{eqnarray} If $A$ is a
matrix representing $g_1$, since $A\in U(2,1)$,
$$AJ_0A^*=J_0.$$ From the fact that $g_1$ fixes $p_1$ and sends
$q_1$ to $p_2$, it is easy to show that $A$ is of the form
\[\left[
\begin{matrix}
   a &0 &0 \\
 a & b & 0 \\
  a(-1-t)/2 &-b& \bar a^{-1}
    \end{matrix}\right]\]
where $|b|=1$.  The fact that $g_1$ sends $q_2$  to $q_3$ implies
that
$$w=bz\bar a-b\bar a=(z-1)b\bar a.$$
$$\frac{-|w|^2+r}{2}=\frac{-|z|^2|a|^2+s|a|^2}{2}+z|a|^2+\frac{-|a|^2-t|a|^2}{2}$$
$$=\frac{-|z-1|^2|a|^2+s|a|^2+z|a|^2-\bar z|a|^2-t|a|^2}{2}.$$
From the first equation we have $a=\frac{b\bar w}{\bar z-1}$ so
$|a|^2=\frac{|w|^2}{|z-1|^2}$. Substituting this into the second
equation gives
$$\frac{-|w|^2+r}{2}=\frac{-|w|^2+(s+z-\bar
z-t)\frac{|w|^2}{|z-1|^2}}{2}.$$ But  this equation is just the
angular invariant identity (\ref{ang}).

The same is true for $g_2$. The equation coming from $g_3$ follows
from the equations from $g_1$ and $g_2$ since three angular
invariants determine the rest in a tetrahedron. In conclusion, out
of 7 parameters $t,s,r,z,w$, there are two angular invariant
equations, which makes 5 dimensional space as expected.

Now we consider holonomy relations. Since there is a unique element
sending three points to other three points in general position
according to the angular invariant, $g_1,g_2,g_3$ are uniquely
determined in terms of $t,z,s,w,r$. Then the holonomy
$g_1^{-1}g_3g_2^{-1}g_1g_3^{-1}\in SU(2,1)$ around  the edge
connecting $\infty$ and $0$ is of the form
\[\left[
\begin{matrix}
   q\alpha &0 &0 \\
 0 &\beta & 0 \\
   0&0& \alpha/q
    \end{matrix}\right]\] where $\alpha$ is a unit complex number
    and $\beta=\alpha^{-2}$.
Now $q\alpha$ is a function in $t,z,s,w,r$. To get a representation
one should have
\begin{eqnarray}\label{ho}
q\alpha=f(t,z,s,w,r)+g(t,z,s,w,r)i=1.\end{eqnarray} Since
$q\alpha=f(t,z,s,w,r)+g(t,z,s,w,r)i$ is a complex number, it adds
two more real equations.  So there are 5 parameters with 2
equations, which gives at most 3 dimensional solution space.   Since
the other holonomy equation for the second edge follows from the
first by Lemma \ref{pure}, the solution space is 3 dimensional.

\end{pf}
Next we deal with $PSp(2,1)$. Here we would like to show that there
is no deformation of representations from the fundamental group of
the figure eight knot complement into $PU(2,1)$, to the ones into
$PSp(2,1)$ out of $PU(2,1)$ up to conjugacy. We do this by
calculating the dimension of the variety of representations near
$\rho_0$ in $PSp(2,1)$ is also 3.

First we begin with a lemma.
\begin{lemma}\label{nonconjugate}
Let $\rho_1,\rho_2:\Gamma \ra SU(2,1)$ be two non-conjugate Zariski
dense representations.  Then they are not conjugate even in
$Sp(2,1)$.
\end{lemma}
\begin{pf}Let $\rho_1(\alpha)=A_1,\rho_1(\beta)=A_2$ and
$\rho_2(\alpha)=B_1,\rho_2(\beta)=B_2$ for the generators $\alpha$
and $\beta$ of the fundamental group of the figure eight knot
complement such that $A_i,B_i\in SU(2,1)$. Suppose $\rho_1$ and
$\rho_2$ are conjugate in $Sp(2,1)$, i.e., there exist $X,Y\in
Sp(2,1)$ such that
$$X A_1 X^{-1}=B_1,\ X A_2 X^{-1}=B_2$$ where $X=X_1+ X_2 j,\
X^{-1}=X_3+ X_4 j$ and $X_1,X_2,X_3,X_4$ are $3\times 3$ complex
matrices. From $XX^{-1}=I$ there is a relation
$$X_1X_3- X_2\bar X_4=I,\ X_1X_4+ X_2\bar X_3=0.$$
Also conjugation relation gives
\begin{eqnarray}\label{1}
X_1A_1X_3- X_2\bar A_1 \bar X_4=B_1,\  X_1A_1X_4+ X_2\bar A_1\bar
X_3=0,\end{eqnarray}\begin{eqnarray}\label{2} X_1A_2X_3 - X_2\bar
A_2\bar X_4=B_2,\ X_1A_2X_4+ X_2\bar A_2\bar X_3=0.\end{eqnarray}
Since $\rho_1$ and $\rho_2$ are not conjugate in $SU(2,1)$, $X_2\neq
0$.

If $X_1=0$, $X=X_2j,\ X_2\in SU(2,1)$. But $$X_2j A_i (X_2 j)^{-1}=
X_2 j A_i (-j)X_2^{-1}=X_2 \bar A_i X_2^{-1}.$$ It is easy to show
that $X_2 \bar A_i X_2^{-1}\notin SU(2,1)$ for generic element
$A_i$. Since $\rho_i$ is Zariski dense, we may assume that $X_2 \bar
A_i X_2^{-1}\notin SU(2,1)$ by choosing generators $\alpha$ and
$\beta$ properly.

Hence $X_1\neq 0\neq X_2$. Then it is easy to show that $X_3\neq
0\neq X_4$.

Since $\rho_1$ is a Zariski dense representation, we can choose
$A_1$ and $A_2$ arbitrarily independent by choosing a different
generators. Then for $X_1$ and $X_2$, there are at most 18 complex
parameters (indeed from $X_iJ_0X_i^*=J_0$ there are less parameters
than 18), whereas from Equations (\ref{1}) and (\ref{2}) there are
at least $9\times 4$ complex equations. This forces that there are
no solutions.
\end{pf}
\begin{Prop}The dimension of the component of the representation
variety of representations from the fundamental group of the figure
8 knot complement to $PSp(2,1)$, containing $\rho_0$, which cannot
be conjugate into $PSp(1,1)$,  is 3 up to conjugacy.
\end{Prop}
\begin{pf}We claim that there are 10 parameters to
choose two tetrahedra. Put one tetrahedron in a standard position
$$p_1=\infty,p_2=0,q_1=(1,ia),q_2=(z,t).$$
There is one degree of freedom for $q_1$ up to $PSp(2,1)$
(corresponding to Cartan angular invariant, or more concretely
corresponding to $a\in\mathbb{R}$). But there is one parameter
family of isometries fixing $p_1,p_2,q_1$ by Lemma \ref{fix}. So
there are 6 degrees of freedom for $q_2$ (since the dimension of the
boundary of $H^2_\bh$ is 7), which makes 7 degrees of freedom to
choose the first tetrahedron.

Once the first tetrahedron is chosen, three vertices of the second
tetrahedron are determined  by its angular invariant according to
the gluing pattern. To choose the last vertex for the second
tetrahedron, it is connected to three vertices to form three
different faces, so their angular invariants are already determined
in the first tetrahedron by gluing pattern. Hence there are $6-3$
degrees of freedom to choose the last vertex for the second
tetrahedron. In conclusion there are total $7+3$ degrees of freedom
to choose two tetrahedra. Note that our parameters are written in
terms of $(1,ia), (z,t)$ for the first tetrahedron, and
$(1,ib),(w,s)$ for the second, so the parameters are in 10
dimensional subspace of $\br\times\br\times \bh\times\bh\times \Im
\bh \times\Im \bh$.

By the previous Lemma \ref{pure}, a holonomy $H_1$ around an edge
being identity gives 7 equations. So we have 10 variables with 7
real equations. Note that $H_1$ is the product of elements in
$PSp(2,1)$ which can be written in terms of $(1,ia),
(z,t),(1,ib),(w,s)$ but $H_1$ can be written as above in terms of
$r\in \br^+,\mu,\nu\in Sp(1)\subset \bh$. So $H_1=id$ produces 7
independent equations in terms of $(1,ia), (z,t),(1,ib),(w,s)$.
%So
%it suffices to show that there are at least 8 independent equations.
%As in the proof of Proposition \ref{complex}, we obtain the same
%equations (\ref{real}) and (\ref{holonomyeqn}) with the complex
%numbers replaced by quaternionic numbers. If $H_1=id$ is the only
%equation to obtain a representation into $Sp(2,1)$, then all the
%representations into $PU(2,1)$ should be obtained from $H_1=id$ also
%by restricting the equations to $\bc$, which was not the case as we
%saw in the proof of Proposition \ref{complex}. So there should be an
%extra equation to obtain a representation into $Sp(2,1)$ as in
%complex hyperbolic case. This shows that there are 7 equations from
%$H_1=id$ and at least one more from the other equation, so total
%there are at least 8 independent real equations.
Since the dimension of the variety into $PU(2,1)$ is already 3 and
two non-conjugate Zariski dense representations in $PU(2,1)$ cannot
be conjugate by an element in $PSp(2,1)$ by Lemma
\ref{nonconjugate}, the dimension of the variety into $PSp(2,1)$
should be also 3. This shows that every representation in $PSp(2,1)$
around $\rho_0$ is conjugate into $PU(2,1)$. So there is no
deformation.
\end{pf}
 We suspect that this component containing a discrete representation
 $\rho_0$
with purely parabolic holonomy for a peripheral group in $PU(2,1)$
is disjoint from the component containing the holonomy
representation of the complete real hyperbolic structure in
$PSp(1,1)$ in the representation variety in $PSp(2,1)$.

\section{Parameters for the character variety in $PU(2,1)$ near
$\rho_0$} We showed that the dimension of the character variety from
the fundamental group $\Gamma$ of the figure eight knot complement
to $PU(2,1)$ near $[\rho_0]$ is 3. In this section we parameterize
this space using angular invariants. We use the notations of
Proposition \ref{complex}. To parameterize the two ideal tetrahedra,
we used five points
$$p_1=\infty,p_2=0, q_1=(1,t),q_2=(z,s),q_3=(w,r)$$ where $z,w\in\bc, t,s,r\in
\im \bc$. From $\mathbb A(q_2,q_1,p_1)=\mathbb A (q_3,p_2,p_1)$, we
had $$\label{ang}\frac{r}{|w|^2}=\frac{s+z-\bar z-t}{|z-1|^2}.$$ A
direct calculation shows that from $\mathbb A(q_1,p_2,p_1)=\mathbb
A(q_2,p_2,q_3)$ we have
$$\text {arg}(\frac{1-t}{2})=\text {arg}(\frac{(|z|^2-s)(|w|^2+r)(|w-z|^2-r-w\bar z+z\bar w+s)}{8}).$$
Hence there were 5 independent parameters out of $t,s,z,w,r$ to
parameterize two ideal tetrahedra according to the gluing pattern.
Finally the holonomy equation around the edge $(p_1,p_2)$ gave two
more real equations, hence the real dimension of the character
variety around $[\rho_0]$ is 3. Here we give these 3 parameters in
terms of Cartan angular invariants.
\begin{Prop}The character variety $\chi(\Gamma,PU(2,1))$ around
$[\rho_0]$ is parameterized by three angular invariants $\mathbb
A(p_1,p_2,q_j),\ j=1,2,3$.
\end{Prop}
\begin{pf}Since $\mathbb A(q_3,p_2,p_1)=\mathbb A(q_2,q_1,p_1)$,
knowing three angular invariants $\mathbb A(p_1,p_2,q_j),\ j=1,2,3$
will completely determine the angular invariants of the first
tetrahedron by cocycle relation (\ref{cocycle}). The angular
invariants $\mathbb A(p_1,p_2,q_j),\ j=1,2,3$ are functions of only
$t,z,s$. Gluing maps $g_1,g_2,g_3$ relates variable $t,z,s$ to the
variable $w,r$. Hence the holonomy map
$g_1^{-1}g_3g_2^{-1}g_1g_3^{-1}$ relates $t,z,s$ variable to $w,r$.
In other words, holonomy map does not create any  relation among
$t,z,s$. This shows that three angular invariants $\mathbb
A(p_1,p_2,q_j),\ j=1,2,3$ are independent. Hence these three
parameters are parametrization of the character variety around
$[\rho_0]$.
\end{pf}
In \cite{Falbel}, it is shown that the coordinates of tetrahedra
corresponding to $\rho_0$ are
$$p_1=\infty, p_2=0, q_1=(1,\sqrt 3),q_4=(-\frac{-1-i\sqrt 3}{2},
\sqrt 3), q_3=(\frac{-1+i\sqrt 3}{2},\sqrt 3)$$ and
$$\mathbb A(p_1,p_2,q_j)=\frac{\pi}{3},\ j=1,2,3.$$
Hence in our coordinates
$[\rho_0]=(\frac{\pi}{3},\frac{\pi}{3},\frac{\pi}{3})$.

\vskip1cm

\noindent     Inkang Kim\\
     School of Mathematics\\
     KIAS\\ Hoegiro 85, Dongdaemun-gu\\
     Seoul, 130-722, Korea\\
     \texttt{inkang\char`\@ kias.re.kr}

\smallskip

     \end{document}